\documentclass[conference]{IEEEtran}

        \usepackage{caption}
  	\usepackage[pdftex]{graphicx}
  	\graphicspath{{../pdf/}{../jpeg/}}
	\DeclareGraphicsExtensions{.pdf,.jpeg,.png}
        \usepackage{multirow}
        \usepackage{multicol}
	\usepackage[cmex10]{amsmath}
        \usepackage{CJKutf8}
	\usepackage{mathabx}
	\usepackage{algorithmic}
	\usepackage{array}
	\usepackage{mdwmath}
	\usepackage{mdwtab}
	\usepackage{eqparbox}
	\usepackage{url}

        \usepackage[ruled,linesnumbered]{algorithm2e}

\begin{document}
\begin{CJK}{UTF8}{gbsn}

\title{\LARGE Explainable MST-ECoGNet Decode Visual Information from ECoG Signal}

 \author{\authorblockN{Changqing JI\authorrefmark{1}, Keisuke KAWASAKI \authorrefmark{3}, Isao HASEGAWA\authorrefmark{3}, Takayuki OKATANI\authorrefmark{1}\authorrefmark{2}}
 \authorblockA{\authorrefmark{1}∗Graduate School of Information Sciences, Tohoku University, Miyagi, Japan 
 \\Email: \{ji.changqing, okatani\}@vision.is.tohoku.ac.jp}
 \authorblockA{\authorrefmark{3}∗Graduate School of Medical and Dental Sciences, Niigata University, Niigata, Japan
 \\Email: \{kkawasaki,isaohasegawa\}@med.niigata-u.ac.jp}

 \authorblockA{\authorrefmark{2}Center for Advanced Intelligence Project, Riken, Tokyo, Japan}}

\maketitle

\begin{abstract}

In the application of brain-computer interface (BCI), we not only need to accurately decode brain signals, but also need to consider the explainability of the decoding process, which is related to the reliability of the model. In the process of designing a decoder or processing brain signals, we need to explain the discovered phenomena in physical or physiological way. An explainable model not only makes the signal processing process clearer and improves reliability, but also allows us to better understand brain activities and facilitate further exploration of the brain. In this paper, we systematically analyze the multi-classification dataset of visual brain signals ECoG, using a simple and highly explainable method to explore the ways in which ECoG carry visual information, then based on these findings, we propose a model called MST-ECoGNet that combines traditional mathematics and deep learning. The main contributions of this paper are: 1) found that ECoG time-frequency domain information carries visual information, provides important features for visual classification tasks. The mathematical method of MST (Modified S Transform) can effectively extract temporal-frequency domain information; 2) The spatial domain of ECoG signals also carries visual information, the unique spatial features are also important features for classification tasks; 3) The real and imaginary information in the time-frequency domain are complementary. The effective combination of the two is more helpful for classification tasks than using amplitude information alone; 4) Finally, compared with previous work, our model is smaller and has higher performance: for the object MonJ, the model size is reduced to 10.82\% of base model, the accuracy is improved by 6.63\%; for the object MonC, the model size is reduced to 8.78\%, the accuracy is improved by 16.63\%.

\end{abstract}

\IEEEoverridecommandlockouts

\IEEEpeerreviewmaketitle


\begin{figure*}[t] 
\centering
\includegraphics[width=7in]{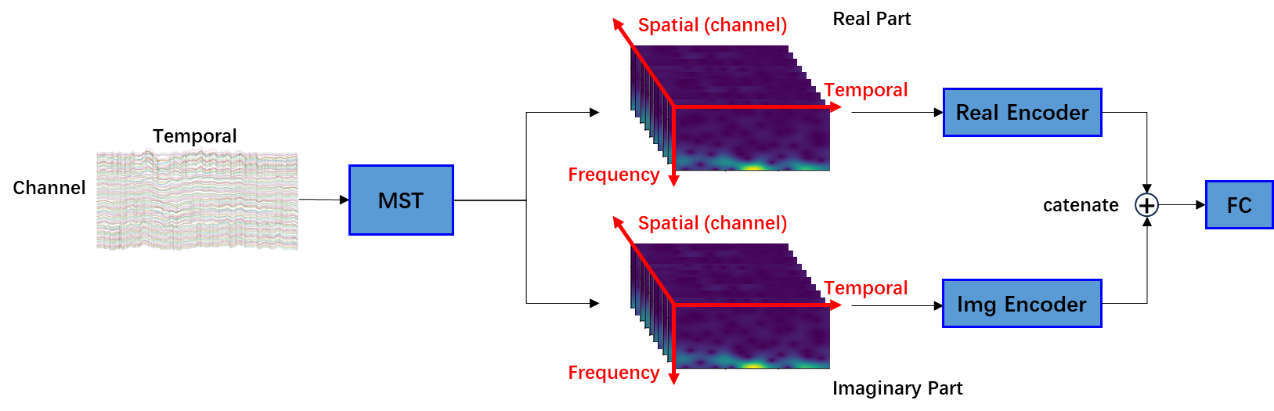}
\caption{\textbf{Outline of Data Flow of MST-ECoGNet}}
\label{fig1:Outline}
\end{figure*}

\section{Introduction}

Humans have never lost the curiosity about brain, the study of brain science cannot be separated from the analysis of electroencephalograms. Electroencephalograms can be divided into two types according to the acquisition method: invasive method, such as ECoG, and non-invasive method like EEG. EEG is a method that attaches the measuring electrodes to the surface of the scalp. It has the advantages of low cost and non-invasiveness. However, the collected data is easily affected by noise, and the electrode layout is sparse in space. ECoG signal acquisition is to place the electrode array directly on the gray matter cortex of the brain, directly collect the electrical signals from the gray matter cortex. Compared with EEG, there is no attenuation of the skull, the signal has a higher signal-to-noise ratio, the high density of electrodes has a high spatial resolution. This makes ECoG widely used, such as epilepsy detection and tracing \cite{2}\cite{3}\cite{4}, brain-computer movement control \cite{1}, and the study of brain visual function \cite{5}\cite{6}\cite{7}. ECoG is a multi-channel waveform temporal series signal. For multi-classification tasks, how to effectively process the data is a challenging work. At the same time, the data processing process also needs to be explainable, which not only helps to improve the reliability of model, but also helps us understand information processing of our brain, provide a certain reference value for future research work.
In the multi-classification task of brain temporal series signals, the model generally can be divided into two parts: feature extractor and classifier.

The feature extractor directly processes the brain signal, extracts the corresponding feature information, including explicit features and implicit features. The explicit features are highly explainable, we can clearly understand the attribute of the features, such as the mathematical statistical information of the brain signal: the average, peak, variance in the temporal domain, and also the frequency, amplitude, phase angle and other information in the frequency domain. \cite{6} used a complex Morlet wavelet to extract temporal-frequency information to study the visual information. \cite{8}\cite{11}\cite{18} used Band power features, \cite{9} used covariance matrices, and \cite{10} used surface Laplacian to complete motor imagery classification. \cite{12} the signal amplitude was used for speller classification. \cite{13}\cite{14} used Power Spectral Density. On the contrary, implicit features always cannot be explained during the extraction process, the feature data exists in the latent space, more like a black box, such as Autoregressive Coefficients (AR) \cite{15}\cite{16}, common spatial patterns (CSP) \cite{17}\cite{19}\cite{20} and feature extraction methods based on deep learning \cite{5}\cite{21}\cite{22}\cite{23}\cite{24}.

The classifier will compress or filter the extracted feature information, then finally classify. General compression or filtering methods include using the peak value or average value, this is common method in the pooling layer of neural networks \cite{23}. There are also mathematical algorithms such as genetic algorithms \cite{25} and principal component analysis (PCA) \cite{26}\cite{27}. The methods for feature classification after filtering include support vector machine (SVM) \cite{30}\cite{31}, linear discriminant analysis (LDA) \cite{32}\cite{33} and k-nearest neighbor (kNN) \cite{28}\cite{29}, decision trees \cite{34}\cite{35} and Naïve Bayes Classifier \cite{36}, as well as neural network methods \cite{37}\cite{38}\cite{39}. Although these methods can complete the final classification task well and efficiently, they always pay less attention on the explainability of data processing.

Based on the above findings and understanding, we focus on designing a model that is explainable and also has high performance. In the feature extraction process of ECoG signals, we use the mathematical method MST \cite{40}. The MST method is an important method for time-frequency analysis of non-stationary signals. Compared with Short-time Fourier transform (STFT) \cite{41}, wavelet transform (WT) \cite{42} and S transform (ST) \cite{43}, it has better resolution in the time-frequency domain both the low-frequency and the high-frequency domain. Furthermore, MST processing of ECoG signals can clearly let us know the physical meaning of the output, its three dimensions represent the frequency domain, temporal domain and spatial domain respectively. For the classifier, we borrowed the idea of EEGNet \cite{23} and adopted a simple neural network. For considering its explainability, we try our best to make the internal module functions have physical or physiological meanings in reality. In this paper, we do not blindly pursue high performance of the model, we pay more attention to make our model to reveal the special information and patterns of ECoG during the process, help us explore the way in which ECoG signals carry visual information. This special information and its intrinsic connections would be an important information for brain scientists to explore and understand the mysteries of the brain. The experiments in this paper reveal that 1) the frequency domain of ECoG signals carries rich information and can provide us the corresponding visual information. Useful information locates in low frequency domain. This finding could have some guiding significance for the design of future sampling frequencies; 2) The Spatial Filter in the classifier can select out effective features, which also means that the visual-based ECoG show some special pattern in spatial domain. This discovery will be helpful for the design of electrode arrays for future ECoG signal acquisition; 3) In the data and result correlation test experiment, we found that there is a 50 msec delay from visual stimuli onset to the appearance of ECoG related information. This discovery may reveal the time-consuming phenomenon of the brain's processing of visual stimuli, which requires further research by brain scientists.

In this paper, introduction of our proposal will be discussed in section II; then model training \& results in section III; in section IV, we present the research experiments; final part is the conclusion.


\begin{table*}[ht]
\centering
\caption{\textbf{Model Structure Details} F means the number of frequency feature map; N means the number of kernel; G means the 'group' parameter in nn.Cov3d function.}
\label{table_structure}
\resizebox{\textwidth}{!}{
\begin{tabular}{l l l l l l}
\hline
Part &Layer & \multicolumn{2}{l}{[Kernel Size] * N / G} &Output &Reference \\
\hline
MST &MST    &\multicolumn{2}{l}{-}  &F x 128 x 300  &ECoG signal size:128 x 300\\
\hline

\multirow{9}{*}{Encoder} &BatchNorm3D &\multicolumn{2}{l}{-} &F x 1 x 128 x 300   \\
{}  &\textbf{nn.Cov3d}&\multicolumn{2}{l}{[1 x 128 x 1] * 2F / F} &2F x 1 x 1 x 300 &\textbf{Spatial Filter}\\ 
{}  &BatchNorm3D &\multicolumn{2}{l}{-} &2F x 1 x 128 x 300   \\
{}  &nn.ELU()    &\multicolumn{2}{l}{-} &2F x 1 x 128 x 300   \\
{}  &nn.AvgPool3d&\multicolumn{2}{l}{[1 x 1 x 4]} &2F x 1 x 128 x 75   \\
{}  &nn.Cov3d    &\multicolumn{2}{l}{[1 x 1 x 16] * 2F / 2F} &2F x 1 x 1 x 75 \\ 
{}  &nn.Cov3d    &\multicolumn{2}{l}{[1 x 1 x 1] * 16 / 1} &16 x 1 x 1 x 75 \\
{}  &BatchNorm3D &\multicolumn{2}{l}{-} &16 x 1 x 1 x 75   \\
{}  &nn.ELU()    &\multicolumn{2}{l}{-} &16 x 1 x 1 x 75   \\
{}  &Flatten     &\multicolumn{2}{l}{-} &1 x 1200   \\
{}  &Cantenate   &\multicolumn{2}{l}{-} &2 x 1200 &Stack real \& img part together   \\

\hline
Full Connection &nn.lin() &\multicolumn{2}{l}{-} &1 x 6

\end{tabular}
}
\end{table*}

\section{Method}

Figure \ref{fig1:Outline} shows the whole step of data processing: the multi-channel temporal series signal ECoG is first analyzed by the MST algorithm to generate a 3D feature space. This feature space is in complex number, which means it has 2 parts: the real part and the imaginary part. These two parts of the data will be processed by the corresponding encoder and then handed over to the full connection (FC) layer for final classification. The encoders of the real and imaginary parts have the same structure, but different learning parameters.

\subsection{Modified S Transform (MST)}

Modified S Transform is an important tool when applying temporal-frequency analysis on non-stationary signals such as ECoG signals. It introduced the arctangent function into the window function, this makes the window width can be adapted to the frequency, ensuring that all parts of the temporal-frequency spectrum have high resolution. This method can effectively improve the energy aggregation with higher time frequency resolution.
MST formula expression as below:\\

\resizebox{0.9\hsize}{!} {$ MST(t, f)= \frac{1}{\sqrt{2\pi}\alpha(f)}\int_{-\infty}^{+\infty}x(\tau)exp(-\frac{(\tau-x)^2}{2\alpha^2(f)})exp(-j2\pi f\tau)d\tau $}
\\
\\
\[
\alpha(f)= \{a[arctan(\frac{f-f_m /2}{b})]+c\}^{-1}
\]

where $f_m$ is the maximum analysis frequency, $f_s$ is the signal sampling rate, parameters $a$, $b$ and $c$ are control factors. Following the parameter selection guideline in MST \cite{40}，in our work we set $ a=5, b=50, c=74, f_m=128, f_s= 1000$.
We use MST to transform ECoG data into 3D feature space. Here we can clearly see that the feature space is composed of temporal domain, frequency domain and spatial domain, which has strong explain-ability.

\subsection{Real \& Imaginary Encoder}

The Real Encoder and the Imaginary Encoder have the same structure. The detailed structure can be found in Table \ref{table_structure}. For high explain-ability, we use very simple structure network. The core part is a single-layer convolution which simulates the Spatial Filter function. In subsequent comparative experiments, we found that the Spatial Filter is the best, which also reveals that the visual-based ECoG signal will show a unique Spatial Pattern at the spatial domain. This pattern carries visual information.

\subsection{Full Connection Layer (FC)}

The full connection layer completes the classification task based on the final features extracted. We concatenate the feature vectors of the Real Encoder and the Imaginary Encoder together, then use the full connection layer to perform six-category classification. Subsequent experimental results show that the parallel structure can effectively utilize the complementary information of the real and imaginary data, can effectively improve the model performance.


\section{Model Training \& Results}
The visual-based ECoG dataset is same dataset used in \cite{5}, we use this dataset to train our proposed model.

\subsection{Visual-Based ECoG Dataset}
\begin{figure}[h!t] 
\centering
\includegraphics[width=3.4in]{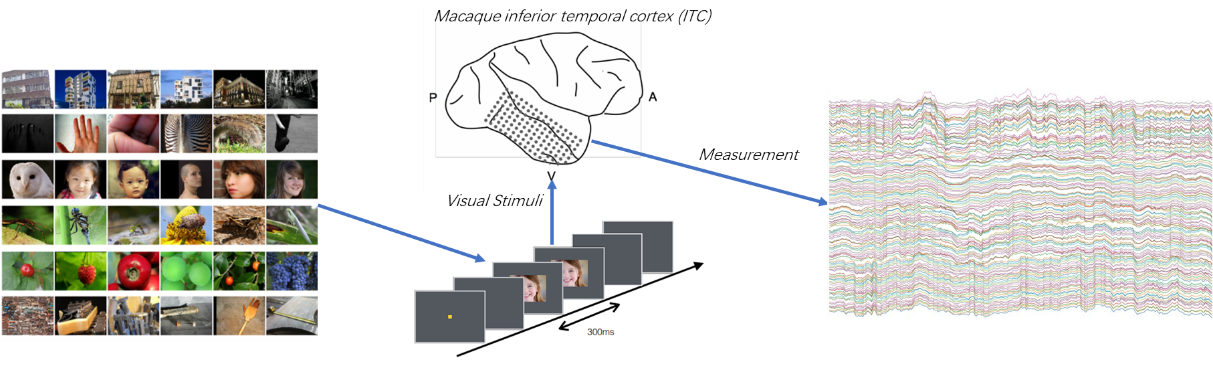}
\caption{\textbf{Outline of ECoG Record} Each image will be used as visual stimuli, last 300 ms, ECoG signal will be measured via electrode array.}
\label{fig2:recoding_outline}
\end{figure}

The ECoG data used in experiment was acquired by Kawasaki et. from Graduate School of Medicine and Dental Sciences, Niigata University. The targets were two macaque monkeys (Subject MonC: 6.1 kg, Subject MonJ: 5.1 kg). The outline of the acquire processing can be seen in Figure \ref{fig2:recoding_outline}. For more information, please refer to \cite{5}\cite{6}\cite{7}. Here we only give an overview:

\subsubsection{Image set selection}
The images used as visual stimuli were selected from natural images, with a total of 6 categories: building, body part, face, fruit, insect, and tool. After screening by three experts, 1,000 images were finally selected for each category, and the image size was reshaped to 512x512;

\subsubsection{ECoG recording}
All animal experiments followed relevant legal requirements, and the experimental subjects were implanted with 8x16 ECoG electrode arrays on the surface of the inferior temporal cortex of the brain, total 128 electrode channels. In fact, subject 2 had another implanted electrode array with of 64 channels. But here we only utilized 128 channels of data which covered the surface of the ITC. In this way, the data of the two subjects remained consistent.

\subsubsection{Recording Process}
The subjects were trained with a visual fixation task to keep their gazes within ±1.5 degree of the visual angle around the fixation target (diameter: 0.3 degree). Eye movements were captured with an infra-red camera system at a sampling rate of 60 Hz. Stimuli were presented on a 15-inch CRT monitor (NEC, Tokyo, Japan) with a viewing distance of 26 cm. In each trial, after 450 ms of stable fixation, each stimulus was presented for 300 ms, followed by a 600 ms blank interval. Signals were deferentially amplified using an amplifier (Plexon, TX, USA or Tucker Davis Technologies, FL, USA) with high- and low-cutoff filters at 300 Hz and 1.0 Hz, respectively. All subdural electrodes were referenced to a titanium screw that was attached directly to the dura at the vertex area. Recording wad conducted at a sampling rate of 1 kHz per channel.

\subsection{ECoG Record Preprocess}
The ECoG record of one trail contains 3 parts, one of which is an active interval that lasts for 300 milliseconds. In this interval, pictures will be displayed to stimulate the monkey's visual system; the other two are the pre-static interval and the post-static interval, pictures are no longer displayed. Please refer to Figure 4 in the data acquisition part. For data preprocessing, we select the 300-millisecond interval before the active interval as the background, the length of the time interval is consistent with the length of the active interval. We calculate the overall mean and variance of the background interval, then use each data of the active interval to subtract this mean and divide it by the variance for normalization. The formula is as follows:
$$ E`_{ac}(n, t)= \frac{E_{ac}(n, t)-m_s}{\sigma_s}$$

$m_s$ $\sigma_s$ indicate overall mean value and standard deviation value of ECoG background respectively. $E_{ac}(n, t)$ indicates the ECoG data in active interval，which has $n$ channels，$t$ sampling points. $E`_{ac} (n, t)$ indicates the ECoG data after processing，this will be the input of MST.
The purpose of this preprocessing is to exclude the influence of the monkey's state during the experiment at different time periods as much as possible. The background interval reflects the basic state of the monkey, if exclude the influence of this period, the remaining changes will be caused by visual stimulation, which is what we are concerned about.

\subsection{Processing of MST}
After above preprocess, we apply MST to perform temporal-frequency analysis, then obtain 3D feature space (temporal-frequency-spatial). The resolution of the frequency domain is 1Hz, and the time domain and space domain are consistent with the original ECoG signal. The object MonC dataset has 21,586 trails, the object MonJ dataset has 13,445 trails, totaling 35,031 trails. In order to obtain more accurate experimental results, the dataset is evenly divided into 5 parts, use 5 cross fold validation to verify the performance of the model. In this paper, all the performance result is the average of 5 times test.

\subsection{Model Training}

The experimental model is implemented in pytorch. The network structure can be found in the encoder part of Figure I and Table I. The data of the two experimental objects were trained separately. The parameters of Real Encoder and Imaginary Encoder were trained on their respective Real and Imaginary datasets. We trained each model for 200 epochs with a batch size of 128. We used the Adam optimizer for optimizing model parameters, with a learning rate of 1.0e-6, a weight decay of 0.0001, β1 = 0.9, and β2 = 0.999. The operating platform information is as follows: cpu Intel Xeon w5-2455x, memory 64 GB, gpu NVIDIA GeForce RTX 4090 24GB.

\subsection{Results \& Comparison}

Table \ref{table_results} shows the experimental test results. Compared with the Base Model \cite{5}, our model performance has been greatly improved: the accuracy of object MonC reached to 53.43\% with an increase of 16.63\%; the object MonJ also increased by 6.63\%. In addition, in terms of model size, our method model is smaller: the experimental object MonC is only 0.0396M, down to 8.78\% of the model size of \cite{5}; the experimental object MonJ is only 0.0488M, down to 1.08\% of the model size of \cite{5}. This shows that our model is smaller, has higher performance, this means more suitable for the application of BCI which requires lightweight models.\\
We analyzed and compared the models of \cite{5}. The reasons for the bloated model, low performance of \cite{5} would be: 1) The feature extractor TCN of \cite{5} contains a complex residual block, processes each channel independently without sharing weight parameters, which makes the number of TCNs large and cannot utilize the correlation between channels; 2) Across-channel transform is applied to 128-dimensional channel data, and 3 matrices (query (Q), key (K), and value (V)) are required for each channel, which makes the parameters of this module too large and is not conducive to parameter learning.

\begin{table}[ht]
\centering
\caption{\textbf{Model Performance Comparison }compare with base work \cite{5}, our model get better performance and smaller size}
\label{table1}
\begin{tabular}{|c|c|c|c|c|c|}
\hline
\multicolumn{2}{|c|}{\multirow{2}{*}{Model}} &\multicolumn{2}{c|}{Params / M} & \multicolumn{2}{c|}{Accuracy}\\
\cline{3-6}

\multicolumn{2}{|c|}{} & MonC & MonJ &  MonC & MonJ  \\
\hline

\multicolumn{2}{|c|}{Base Model \cite{5}} &0.451 &0.451	&36.8±0.51 &27.59±0.73\\ 
\hline

\multicolumn{2}{|c|}{\textbf{MST-ECoGNet}}  &\textbf{0.0396} &\textbf{0.0488} &\textbf{54.15±1.15} &\textbf{35.98±1.15}\\ 
\hline

\end{tabular}
\label{table_results}
\end{table}

\section{Research Experiments}
There are few studies on ECoG signals based on the visual system. Thanks to the excellent work of the biology team at Niigata University, we are fortunate to have access to these unique data. In order to better and more deeply understand the ECoG signal, here we design various experiments to explore what features in the ECoG signal can carry visual information. We studied the three dimensions (temporal-frequency-spatial) of the 3D feature space and also the two data types (real part and imaginary part). These findings are also the basis for the designing of our proposed model.

\subsection{Influence of Frequency Range}
\begin{figure}[h!t] 
\centering
\includegraphics[width=3.4in]{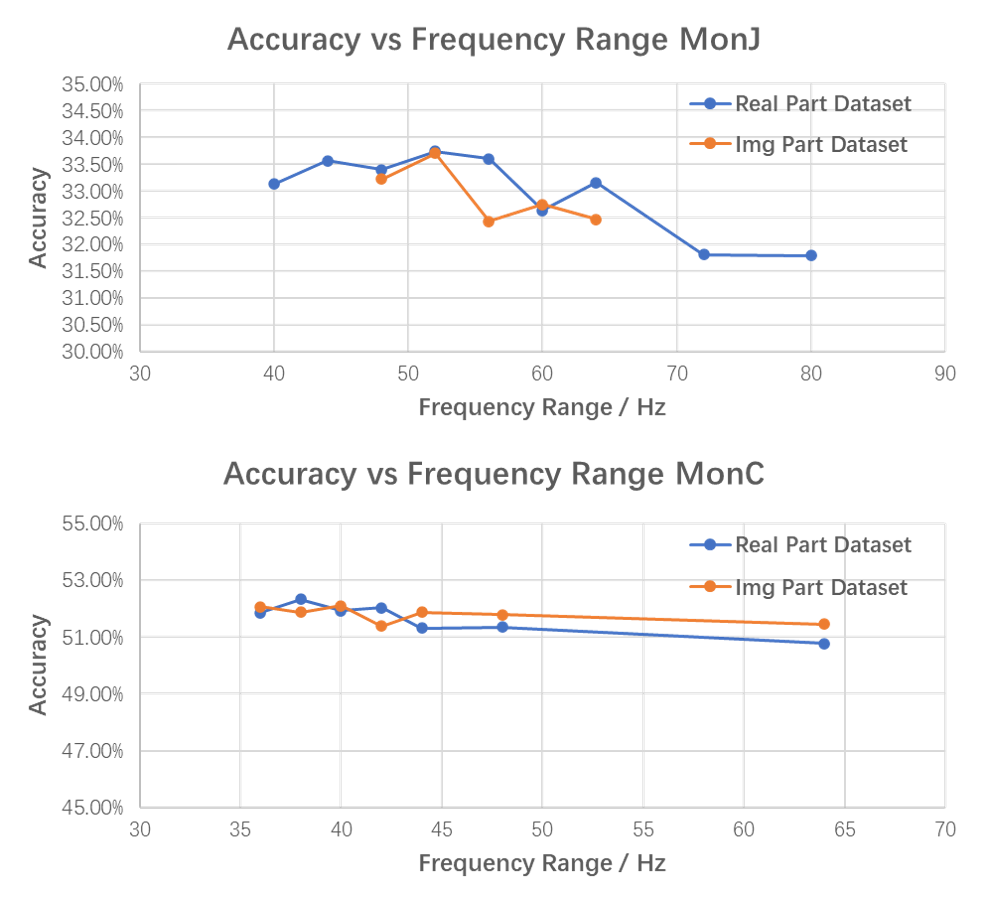}
\caption{\textbf{Influence from Frequency Range} Results test on 2 subjects \& 2 data part, final frequency range: MonJ 0-52Hz, MonC 0-38Hz.}
\label{fig3:fig3_frequencyRange}
\end{figure}

In order to better obtain the features in frequency domain, when MST processes ECoG signals, we set the resolution of the frequency domain to 1Hz, and the frequency range is up to 128Hz. However, we observed that the signal strength on the temporal-frequency diagram is mainly concentrated in the low-frequency part, which requires us to reduce the frequency range. This has two advantages: reducing the frequency range will make the effective data more concentrated, which is conducive to the model to better extract features and improve performance; reducing the frequency range will make the data set smaller, which is conducive to improving the efficiency of model training. To this end, we designed a frequency range screening experiment: we reduce the frequency range step by step, each step we use the reduced dataset to train the network. The network here used is consistent with that shown in Figure I, but only one encoder is used, and the input dimension of the corresponding FC is also reduced by half. We decide the best frequency range based on the performance of the network. In the experiment, two experimental subjects were tested respectively, we also separately test the 2 types of the data: Real part data and Imaginary part data. The experimental results are shown in Figure \ref{fig3:fig3_frequencyRange}.

The experimental results clearly show the optimal value of the frequency range. The final range of object MonJ is 0-52Hz, the frequency range of object MonC is 0-38Hz. At the same time, it also confirms that effective information is contained in the low-frequency range, frequency is also one of the important parameters affecting model performance.

\subsection{Frequency Temporal Spatial Feature}

\begin{figure}[h!t] 
\centering
\includegraphics[width=3.4in]{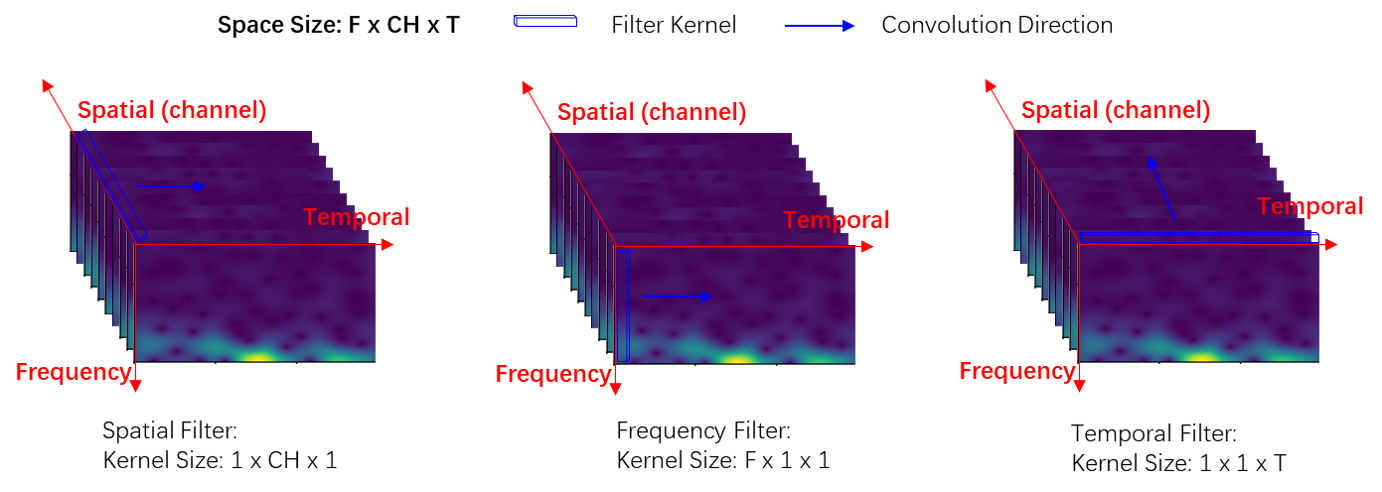}
\caption{\textbf{Filter Structure for 3 Different Dimensions} Each filter has its own kernel size and convolution direction, then the filter can only focus its own dimension.}
\label{fig4:filtersStructure}
\end{figure}

After the MST temporal-frequency analysis of the ECoG signal, we accurately obtain the feature information of three dimensions. We want to know whether information of these three dimensions carry visual information. To this end, we made some modifications to the structure of the encoder, designed three different filters, each corresponds to different dimensional features. The structure of each filter is shown in Figure \ref{fig4:filtersStructure}. Each filter has its own focused dimension, through comparing the respective performance results, we can clearly know which dimensional information carries more visual information.

\begin{figure}[h!t] 
\centering
\includegraphics[width=3in]{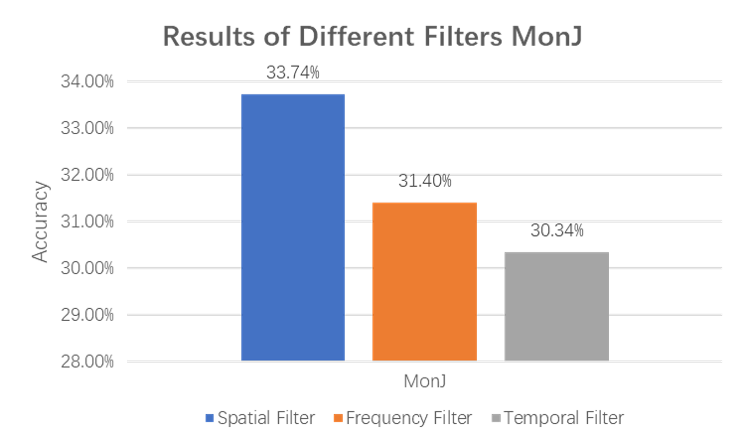}
\caption{\textbf{Results of 3 Different Filters} We only test on MonJ dataset, real part. Results show Spatial dimension has best performance, which means spatial pattern carries visual information.}
\label{fig5:filtersResults}
\end{figure}

Experimental results refer to Figure \ref{fig5:filtersResults}. The results show that three different dimensional filters can extract effective feature information to complete the classification task, but the performance of using the Spatial Filter is better. This shows that there is a Spatial Pattern in the ECoG signal, which contains visual information. This discovery is worth sharing with the biology group at Niigata University, this would be of reference significance in the design shape of the electrodes array or the resolution between electrodes in the future.

\subsection{Amplitude-Angle vs Real-Imaginary}
\begin{figure}[h!t] 
\centering
\includegraphics[width=3in]{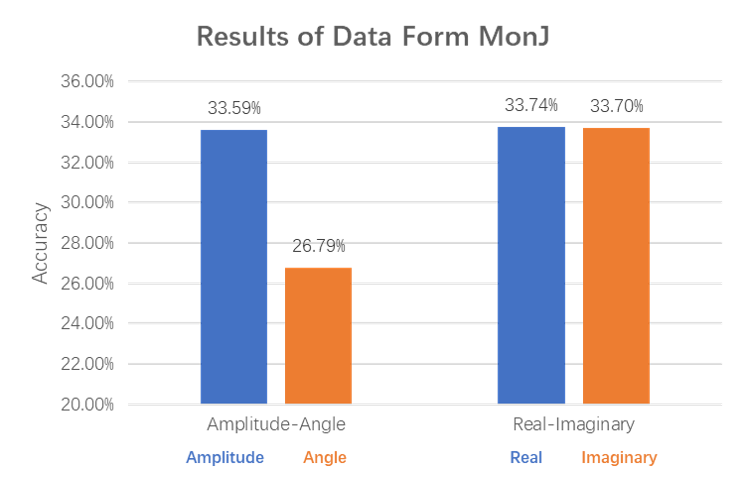}
\caption{\textbf{Results of 2 Different Data Form} Only test on MonJ dataset, Results show data form of Real-Imaginary would be better than Amplitude-Angle data form.}
\label{fig6:fig6_DataFormResults}
\end{figure}

The output of MST processes is complex data. For complex data, we can use them in two forms: Amplitude-Angle \& Real-Imaginary data. We tested these two forms of data sets respectively; the results of the comparative experiment are shown in Figure \ref{fig6:fig6_DataFormResults}. We can clearly see that compared with the data type of amplitude-angle, the real-imaginary data type has better performance.

\begin{figure}[h!t] 
\centering
\includegraphics[width=3.5in]{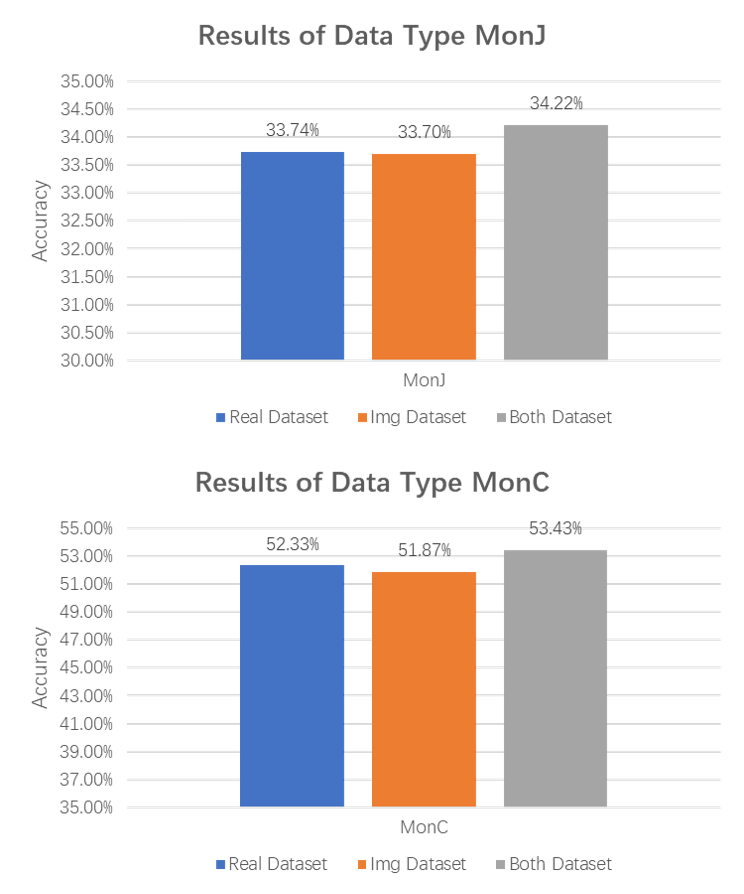}
\caption{\textbf{Parallel Encoder vs Single Encoder} Parallel encoder using both real and imaginary data has better performance compare with single encoder which only one data real or imaginary part.}
\label{fig7:fig7_ResultsDataType}
\end{figure}

In addition, considering the Euler formula, the two forms of data are equivalent, the amount of information contained should be same. However, the Real-Imaginary data type is more conducive for the model to extract effective information.
$$ e^{-j\theta}=cos⁡(\theta)+jsin(\theta)$$
The real and imaginary parts of the complex domain are independent of each other, the feature space composed of real and imaginary is also independent with each other, so the information contained in these two parts should also be complementary. Use a parallel structure to process real and imaginary data at the same time, concatenate the extracted feature vectors, then classify them via the FC layer. This can effectively utilize the complementary information of the two parts of the data. To this end, we designed a comparative test: single encoder and parallel encoder. The experimental results are shown in Figure \ref{fig7:fig7_ResultsDataType}:\\

The comparison of experimental results shows that using both data in parallel can effectively improve the performance. In addition, in the frequency importance test, we also found that the real and imaginary data show a complementary trend in the low frequency domain.

\subsection{Timing for Distinguish Information}
How long does it last from the start of visual stimuli onset to the point of ECoG starts to carry visual information? This is an interesting topic among the experts of the biology group. In order to obtain the corresponding data, we designed a test experiment: first, we used the 3D Feature data obtained by the MST method to train the corresponding model, then we modified the data set, retaining the data of 5 sequent samples period along the temporal axis, and setting all the other period data to zero. The modified data set was tested for accuracy on the pre-trained model, the accuracy obtained was used to indicate the importance of the data in this time period. The higher the importance, the more visual information this time period carries. For detailed experimental algorithms, please refer to Algorithm \ref{alg:alg_tmpfrq}.\\

\begin{figure}[h!t] 
\centering
\includegraphics[width=3.4in]{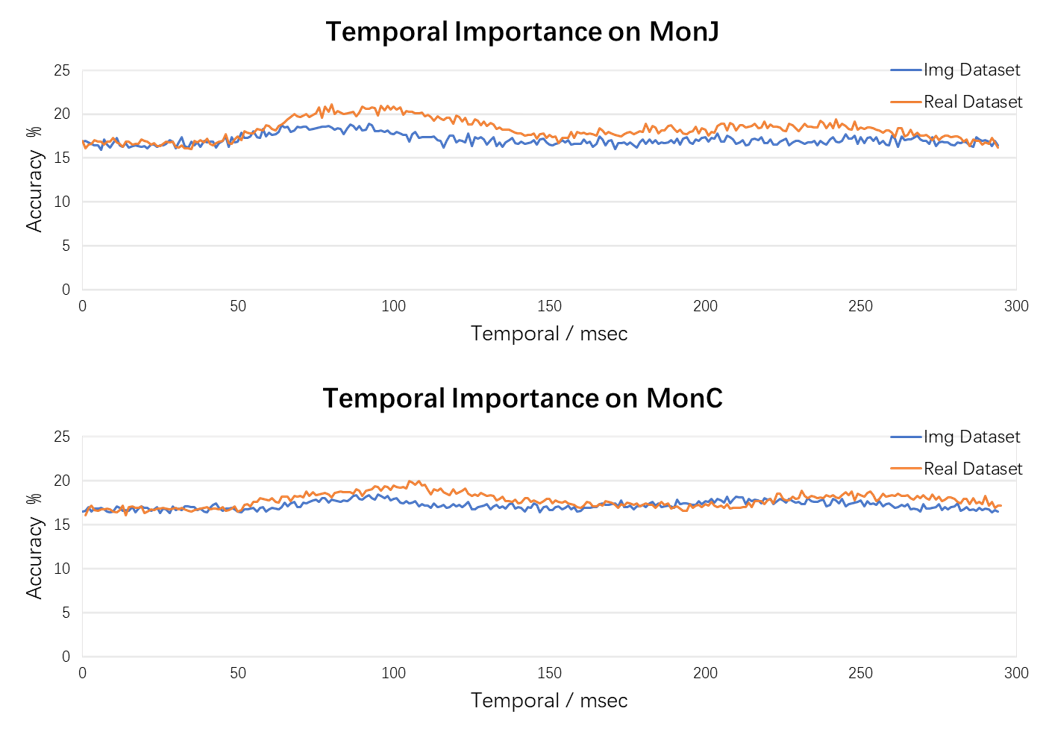}
\caption{\textbf{Temporal Importance Test} Results show that along temporal axis, there is clear time lag about 50 msec in both test subjects.}
\label{fig8:fig8_TemporalImportance}
\end{figure}

Figure \ref{fig8:fig8_TemporalImportance} shows the importance of time domain data. We tested the data sets of two different objects separately, and also tested the real and imaginary data separately. As can be seen from the data in the figure, ECoG starts to carry visual information about 50 milliseconds after the visual stimuli onset, and it shows two waves in the trend. Similar phenomena occurred in both experimental subjects. The mechanism of this phenomenon needs to be explained with the help of the professional knowledge of the biology group. We will share the phenomena we found with the biology group.

\subsection{Complementary Information between Real \& Imaginary Data}
\begin{figure}[h!t] 
\centering
\includegraphics[width=3.4in]{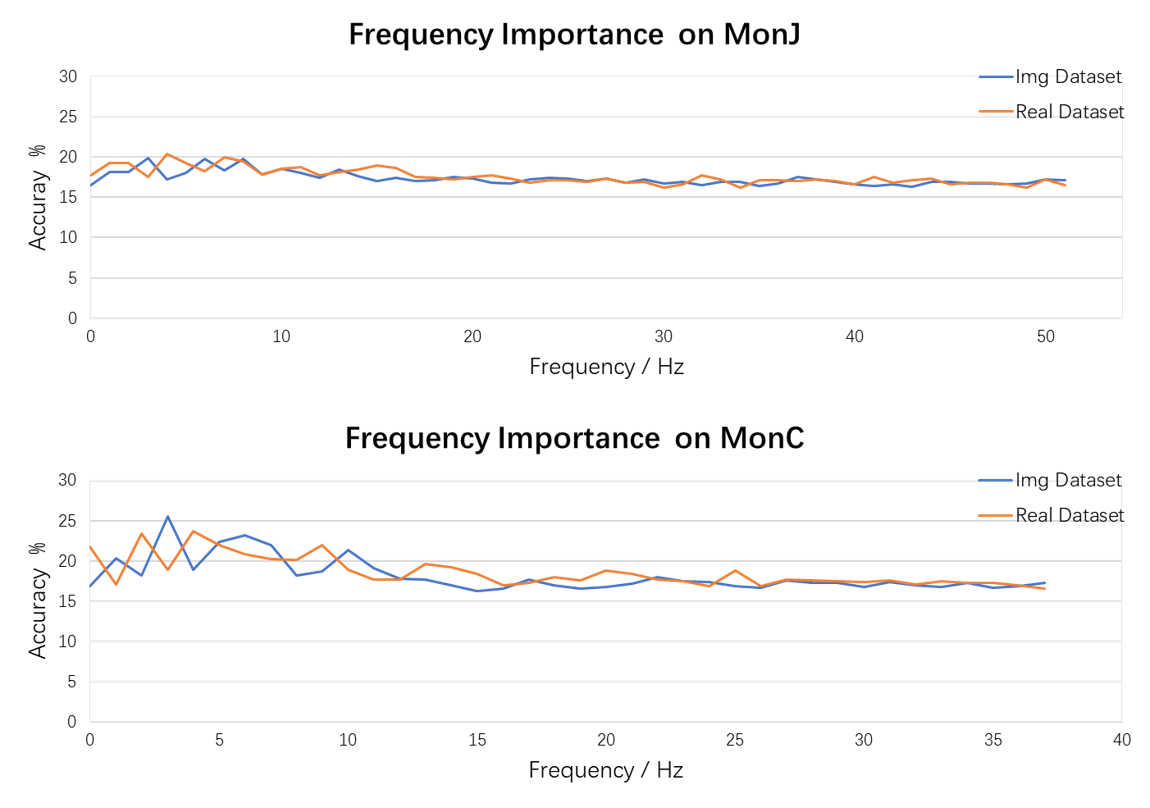}
\caption{\textbf{Frequency Importance Test} Results show that most information comes from low-frequency region, the complementary phenomenon exists between real and imaginary data.}
\label{fig9:fig9_FrequencyImportance}
\end{figure}

The parallel processing of real and imaginary data improves the performance of the model, theoretically the feature information provided by these two parts of data will be independent and complementary. In this regard, we designed a verification experiment, the method is the same as section 4, the only difference is that for each step only one frequency of feature data is retained. Similarly, two objects and two data parts were tested. The specific experimental results are shown in Figure \ref{fig9:fig9_FrequencyImportance}.\\

The experimental results clearly show that in the low-frequency domain, the real and imaginary data show complementary characteristics, which confirms our conjecture and proves that the performance of the model we proposed is reliable.

\begin{algorithm}[h!t]
    \caption{Temporal (Frequency) Importance Test}
    \label{alg:alg_tmpfrq}
    \KwIn{original ECoG dataset $E_{ori}$}
    \KwOut{Model accuracy list $Acc$=$[acc_0,…,acc_i,…,acc_n]$,  $acc_i$ means model accuracy after only keep $i^{th}$ temporal segment based on original 3D feature space}
    \# each temporal segment has 5 time points\\
    initialization\\
    $Acc$=$[ ]$, temporal start point $t_i=0$, end time point $t_{off}= t_i+5$  \\
    Pre-trained model $Model$
    
    \While{$t_i < $ 295}{
        
        for Dataset $E_{ori}$, set temporal period $[0: t_i] \& [t_{off} : 300]$ to 0, keep rest part, then get new ECoG dataset $E^*_i$\\
	$acc_i=Model(E^*_i)$   			\\
 	$Acc.append(acc_i)$				\\
	$t_i +=1$			
       }
    \# for frequency importance test, use same step, but each step we only keep one frequency feature map.
\end{algorithm}


\section{Conclusion}
There are very few studies on ECoG signals based on visual stimuli. Thanks to the excellent work of the biology team at Niigata University, we are fortunate to have access to these precious data. As a hot topic: Is there a corresponding connection between visual stimuli and ECoG signals? Here we use the ECoG signal multi-classification task to start the study. Taking this as an opportunity, we further designed experiments to explore what kind of features carry visual information in multi-channel temporal series signals such as ECoG.\\
MST is a good temporal-frequency analysis method, which provides us a way of three different dimensions to observe and analyze ECoG signals: temporal domain, frequency domain, and spatial domain. Compared with the original ECoG signal, which only has spatial domain and temporal domain, there is an additional frequency domain. From the signal intensity map generated by MST, we can see that the main information is concentrated in the low-frequency band, and the best frequency band needs to be solved first. This is related to the performance of the model and the timeliness of model training. The results of the experiment verified that the low-frequency band is more effective, and it was also found that the best frequency band has individual differences, which also requires further brain research.\\
In addition, to further study the 3D feature space, we designed an extremely simple and explanatory network structure to process the 3D feature space along different dimensions. The core of this network is a network composed of a single layer of convolution in the encoder, which we call Filter. We designed three different convolution kernels with corresponding specific convolution directions, which correspond to the three dimensions of the 3D feature space. From the test results, we can see there is a great improvement compared with \cite{5}. Among them, the result of Spatial Filter is the best. This result reveals that the frequency domain as well as the spatial domain carry the visual information. At the same time, there is a specific pattern in the spatial domain of ECoG, and this pattern can be used to extract visual information. This finding suggests that we need to consider the design of the electrode array including shape and density in the future.\\
In addition, the output data of MST is a complex structure, which has two expression forms: Amplitude-Angle and Real-Imaginary. We designed an experiment to compare the impact of these two forms of data on the performance of the model. The experimental results show that the data form of Real-Imaginary is more conducive to the extraction of features by the neural network. In addition, from the Euler formula, we know that the amount of information contained in the two forms of data is consistent, and the data information of Real and Imaginary part are independent with each other. The complementary features that can be extracted by using Real and Imaginary data in parallel, this is conducive to improving model performance. The experimental comparison results verified this conjecture. Compared with using Amplitude data alone, which is generally used in other models, using Real and Imaginary in parallel will be better.\\
Furthermore, for the impact of features in temporal domain, we designed a test experiment of Temporal Importance. The data set only retains a short period of data along the temporal axis, and the rest is set to zero. The processed data is tested on the trained model for accuracy, using this way we can find out which period of the data can provide useful information. That means that period ECoG carries visual information. The experimental results show that both experimental subjects have a lag time of about 50 milliseconds compared to visual stimuli, two waves appear. This discovery is worth sharing with biological experts and maybe inspire future brain research.\\
Finally, in this paper we used a simple method to explore the unique characteristics of ECoG signals from three dimensions: frequency domain, time domain, and spatial domain, and explored its possible ways of carrying visual information. These findings provide new directions for future data processing or research, also provide potential discussion directions for biological experts. At the end of the last, based on the above findings, we proposed an explanatory model MST-ECoGNet, which is lighter and has better performance than the model in \cite{5}.


\section*{Acknowledgment}
The author of this paper is financially support by CSC (China Scholarship Committee) and RA of OKATANI Lab in GSIS of Tohoku University.


\bibliographystyle{IEEEtran}
\bibliography{IEEEabrv,biblio}


\end{CJK}
\end{document}